\documentclass[12pt]{amsart}
\usepackage{amsmath,amssymb,amsthm,amsfonts,youngtab,graphicx}
\usepackage{epsfig}
\usepackage[curve]{xy}
\newtheorem{theorem}{Theorem}[section]
\theoremstyle{definition}
\newtheorem{definition}[theorem]{Definition}
\newtheorem{lemma}[theorem]{Lemma}
\newtheorem{proposition}[theorem]{Proposition}
\newtheorem{corollary}[theorem]{Corollary}
\newtheorem{example}[theorem]{Example}

\theoremstyle{remark}

\numberwithin{equation}{section}

\begin{document}

\title{On the diagonal hooks of a symmetric partition}
\author{Rishi Nath}
\address{York College/City University of New York}
\email{rnath@york.cuny.edu}

\subjclass[2000]{203C0}


\date{}

\dedicatory{}

\begin{abstract}
Using only a symmetric $p$-core and $p$-quotient, we give an
explicit formula for the set of diagonal hook lengths of the
associated symmetric partition.
\end{abstract}
\maketitle
\section{Introduction}
Suppose $\mathbb{N}=\{0,1,\cdots\}$ let $n\in\mathbb{N}$ and $p$ be
a prime. For standard definitions of a {\it partition} $\lambda$ of
$n$, its {\it dual} $\lambda^*$ and {\it Young diagram} $[\lambda]$,
a {\it hook} $h_{ij}$ of $[\lambda]$ with {\it corner} $(i,j)$, the
{\it hook length} $|h_{ij}|$, the {\it arm length} and {\it leg
length} of $h_{ij}$, and a $\beta$-set $X$ corresponding to
$\lambda,$ we refer readers to \cite{F-S}, \cite{J-K},
\cite{Mac}, \cite{Ols}.

A $\beta$-set $X$ associated to a partition $\lambda$ can be seen as
a finite set of non-negative integers, represented by beads at
integral points of the $x$-axis, i.e. a bead at position $x$ for
each $x$ in $X$. Then $X$ is a $\beta$-set to $\lambda$ in the {\it
extended sense} if we extend $X$ infinitely in both directions with
beads at all negative positions and spaces at all positions to the
right of the position of the largest integer $x_k\in X$. In this
interpretation, $\beta$-sets equivalent to $X$ are the same infinite
string of beads and spaces with the origin shifted a finite number
of positions to the left. A {\it minimal} $\beta$-set $X$ is an
extended set where the first space is counted as 0. If $X$ is a
minimal $\beta$-set of $\lambda$, we define $|X|$ as the number of
beads occurring to the right of the leftmost space.

Given a fixed integer $p$, we can arrange the nonnegative integers
in an array of columns and consider the columns as runners of an
{\it abacus} in order to represent $X$.
\[
\begin{array}{cccc}
0      &  1     & \cdots &  p-1\\
p      &  p+1   &  &  2p-1\\
\vdots  &   & \ddots & \\
mp &  & &  mp+p-1
\end{array}
\]
The column containing $\gamma$ for $0\leq{\gamma}\leq{p-1}$ will be
called the {\it $\gamma$th runner} of the abacus. The integers
$\gamma,\gamma+p,\gamma+2p,\cdots$ label corresponding positions
$0,1,2,\cdots$ on the $\gamma$th runner. Placing a bead at position
$x_{j}$ for each $x_{j}\in{X}$ gives the $\textit{abacus diagram}$
of $X$.

We deviate from standard notation by denoting by $h=(y,x]$ a hook
arising from the $\beta$-set $X$ of $\lambda$ where $y\not\in X$ and
$x\in X$. We define the {\it hook length} of $(y,x]$ as $x-y.$ Lemma
\label{bhama} describes a bijection between the set of hooks
$h_{ij}$ of the Young diagram $[\lambda]$ and the set of hooks
$(y,x]$ of a $\beta$-set $X$ of $\lambda$.
\begin{lemma}\label{bhama} Let $\lambda$ be a partition of $n$ and $X$ a $\beta$-set of $\lambda$.
A hook $h=(y,x]$ of $X$ corresponds to the hook $h_{ij}$ with corner
node $(i,j)$ in the Young diagram $[\lambda]$ where
\[
i=|z\in\mathbb{N}:z\in X, z\geq x|
\] and
\[
j=|z\in \mathbb{N}:z\not\in X, z\leq y|.
\]
Additionally, the leg length and arm length of $h$ are $|z\in
\mathbb{N}:z\in X,y<z<x|$ and $|z\in \mathbb{N}:z\not\in X:y<z<x|$
respectively.
\end{lemma}
\begin{proof} See pg. 180 in \cite{F-S}.
\end{proof}

If $h=(y,x]$ is a hook of length $p$ (henceforth a {\it p-hook}) of
$X$ then $\{y\}\cup X-\{x\}$ is a $\beta$-set for a partition
$\lambda_{1}$ of $n-p.$ We say that $\lambda_{1}$ and
$X_{\lambda_{1}}$ are achieved from $\lambda$ and $X$ respectively
by removing a $p$-hook. In the opposite manner, we see that
$\lambda$ and $X$ are gotten from $\lambda_{1}$ and
$X_{\lambda_{1}}$ respectively by adding a $p$-hook. Subsequently,
the abacus diagram of $X_{\lambda_{1}}$ is related to that of $X$ by
moving the bead at $x\in{X}$ up one position on the runner. Let
$X^0$ be the unique $\beta$-set obtained from $X$ by successively
removing $p$-hooks until none are left. Thus $X^0$ will have no
$p$-hooks. The partition $\lambda^0$ represented by $X^0$ is called
the {\it p-core} of $\lambda$ and is uniquely determined by
$\lambda$.  The abacus of the $p$-core $\lambda^0$ is obtained from
the abacus of $\lambda$ by pushing up the beads in each runner as
high up as they can go (Theorem 2.7.16,\cite{J-K}).

A hook $h=(y,x]$ of length divisible by $p$ is said to be on the
$\gamma$th runner if $x$ is on the $\gamma$th runner. Then $y$ is
also on the $\gamma$th runner. In particular, hooks of length
divisible by $p$ are on the same runner if and only if they have the
same residue modulo $p$. For $0\leq{\gamma}\leq{p-1},$ let
$X_{\gamma}=\{j:\gamma+jp\in{X}\}$ and let $\lambda_{\gamma}$ be the
partition represented by the $\beta$-set $X_{\gamma}$. Notice that
this is the partition whose beads appear on the $\gamma$th runner of
the abacus diagram of $\lambda.$ Our convention will be that the
{\it p-quotient} of $\lambda$ is the sequence
$(\lambda_{0},\cdots,\lambda_{p-1})$ obtained from $X$ where
$|X|\equiv 0\pmod{p}$. We call $X_{\gamma}$ the $\beta$-set of
$\lambda_{\gamma}$ {\it induced} by $X$.

A partition is {\it symmetric} if $\lambda=\lambda^*.$  A
$p$-quotient $(\lambda_0,\cdots,\lambda_{p-1})$ is {\it symmetric}
if $\lambda_{i}=\lambda^*_{p-i-1}$ where $0\leq i\leq p-1.$ The
$p$-quotient and $p$-core of a partition $\lambda$ and its dual
$\lambda^*$ are related in the following manner.
\begin{lemma}\label{sym_quo}Let $X$ be a $\beta$-set for $\lambda$  such that
$|X|\equiv{0}\pmod{p}$. Let $\lambda^*$ be the dual of $\lambda$,
let $(\lambda^*)^0$ be the $p$-core of $\lambda^*$ and let
$\{\lambda^*\}_{}=\{\lambda^*_0,\cdots,\lambda^*_{p-1}\}$ be the
$p$-quotient of $\lambda^*$. Then $(\lambda^*)^0=(\lambda^0)^*$ and
$(\lambda_{\gamma})^*=\lambda_{p-1-\gamma}$ for $0\leq\gamma\leq
p-1$. Hence $\lambda=\lambda^{*}$ if and only if
$\lambda^0=(\lambda^0)^*$ and
$(\lambda_{\gamma})^*=\lambda_{\gamma^*},$
\end{lemma}
\begin{proof} See Proposition 3.5 in \cite{Ols}.
\end{proof}
Given a symmetric partition $\lambda$, we let
$\delta(\lambda)=\{\delta_{ii}(\lambda)\}$ be the set of diagonal
hooks where $h_{ii}=\delta_{ii}(\lambda)$. When there is no
ambiguity we will set $\delta_{ii}(\lambda)=\delta_{ii}$. By abuse
of notation $\delta_{ii}$ will also stand for the size of
$\delta_{ii}(\lambda)$, and $\delta(\lambda)$ for the set of
diagonal hook lengths of $\lambda$.

In this paper we give an explicit formula for
$\delta(\lambda)=\{\delta_{ii}\}$ in terms of only the $p$-quotient and
the $p$-core. These results are motivated by an ongoing study of the
irrationalities of the character table of the alternating groups
$A(n)$, which, by a classical result of Frobenius, arise from the
diagonal hook lengths of symmetric partitions of $n$.  In particular
they assist in verifying that a recent refinement of McKay's
conjecture by Navarro \cite{Nav} involving Galois automorphisms
holds for $A(n)$ in special cases \cite{Nat}.
\section{Bisequences and diagonal hooks} A {\it
bisequence}
\[
(\alpha_1,\alpha_2\cdots,\alpha_t|\beta_1,\beta_2\cdots,\beta_t)
\]
will be an ordered pair of strictly decreasing sequences of
non-negative integers
\[
\begin{array}{l}
(\alpha_1,\alpha_2\cdots,\alpha_t)\\
(\beta_1,\beta_2\cdots,\beta_t)
\end{array}
\]
of the same length $t$. For example let $\alpha_i$ and $\beta_i$ be
the leg and arm lengths of $\delta_{ii}\in \delta(\lambda)$. Then
the sequences $(\alpha_1,\alpha_2,\cdots,\alpha_t)$ and
$(\beta_1,\beta_2,\cdots,\beta_t)$, are strictly decreasing. Hence
we may define the bisequence
\[
D(\lambda)=(\alpha_1,\alpha_2,\cdots,\alpha_t|\beta_1,\beta_2,\cdots,\beta_t).
\]
We define the {\it components} of $D(\lambda)$ to be
$D(\lambda)_{L}=(\alpha_1,\alpha_2,\cdots,\alpha_t)$ and
$D(\lambda)_{R}=(\beta_1,\beta_2,\cdots,\beta_t)$. An element of
$D(\lambda)$ is an ordered pair $(\alpha_i|\beta_i)$ for some $i$
and corresponds to a diagonal hook $\delta_{ii}$. Then
$|D(\lambda)|$ is the number of such pairs, and equals $t$. Note
that $|D(\lambda)_R|$ and $|D(\lambda)_L|$ (the number of arm
lengths and leg lengths of the diagonal hooks respectively) both
equal $t$ as well. The {\it dual} of $D(\lambda)$ is
\[
D(\lambda)^*=(\beta_1,\beta_2,\cdots,\beta_t|\alpha_1,\alpha_2,\cdots,\alpha_t).
\]
Clearly $D(\lambda)^*=D(\lambda^*)$. If $\lambda$ is symmetric then
$D(\lambda)_{L}=D(\lambda)_R$. We attach to $D(\lambda)$ a $p$-tuple
$D'(\lambda)$ of bisequences.
\begin{definition}
Let $D'(\lambda)=(D_0(\lambda),\cdots,D_{p-1}(\lambda))$ where
$D_{\gamma}(\lambda)$ is defined as follows.
\begin{enumerate}
\item If $\alpha=\gamma+mp\in D(\lambda)_{L}$,
$0\leq\gamma\leq p-1$ and $m\geq 0$, we put $m$ in
$D_{p-1-\gamma}(\lambda)_{L}$.
\item If $\beta=\gamma+mp\in D(\lambda)_{R}$,
$0\leq\gamma\leq p-1$ and $m\geq 0$, we put $m$ in
$D_{\gamma}(\lambda)_R$.
\end{enumerate}
$D'(\lambda)$ is called the {\it p-quotient} of $D(\lambda)$.
\end{definition}

From this definition, given $D'(\lambda)$, we can obtain
$D(\lambda)$. At the moment, it is not clear that for each $\gamma$
the sequences in $D_{\gamma}(\lambda)_L$ and $D_{\gamma}(\lambda)_R$
have the same length. This will be shown to be true when
$\lambda^0=\emptyset$ in Theorem \ref{eqlgth}.
\section{$\theta(\lambda)$ and diagonal hooks}
The diagonal hooks $\delta_{ii}$ of $\lambda$ correspond to the
following hooks of $\beta$-set $X=\{x_1,\cdots,x_k\}$. The largest
hook $\delta_{11}$ corresponds to $(y_1,x_k]$, where $y_1$ is the
position of the smallest space i.e. the minimal positive integer not
included in $X$. By removing $\delta_{11}$ (that is, by moving the
bead at position $x_k$ to the space $y_1$) then $\delta_{22}$
corresponds to the largest hook of $\lambda^{\vee}$, and so on. Thus
the diagonal hooks correspond to the nested hooks starting with the
longest hook in $X$, then the longest hook contained strictly within
that longest hook, and so on.

Let $\lambda^{\vee}$ be the partition obtained from $\lambda$ by
removing $\delta_{11}$. Then $X^{\vee}$ is the induced $\beta$-set of
$\lambda^{\vee}$.
\begin{proposition}\label{nsymax} Suppose $\lambda$ is a partition of $n$ and
let $X$ be a $\beta$-set for $\lambda$. Then there exists a
half-integer $\theta(\lambda)$ such that the number of beads to the
right of $\theta(\lambda)$ equals the number of spaces to the left
of $\theta(\lambda)$.
\end{proposition}
\begin{proof} Let $\theta$ be the point at the half-integer
just to the left of the smallest space of $X$. So there are 0 spaces
to the left of $\theta$ and a finite number of beads to the right.
Next move $\theta$ a unit distance to the right, that is, to the
next half-integer on the right. One and only one of the following
happens: Either the number of spaces to the left of $\theta$
increases by one or the number of of beads to the right decreases by
one. So by iterating this process we reach a point where the number
of spaces to the left of $\theta$ is equal to the number of beads to
the right of $\theta$. That point is then $\theta(\lambda).$
\end{proof}

Let $X$ be a $\beta$-set for a partition $\lambda$ (not necessarily
symmetric), with maximal element $x_k\in X$. Let $X_+$ be the subset
of beads to the right of $\theta(\lambda)$ and let $X_-$ be the
subset of spaces to the left of $\theta(\lambda)$. We index the
elements of $X_+=\{y'_i:i\leq r\}$ so that $y_r'<\cdots<y_1'$ with
$y_1'$ as the largest bead. Correspondingly, we index the elements
of $X_-=\{y_i:i\leq r\}$ so that $y_1<\cdots<y_r$ with $y_1$ as the
smallest space. In particular, $y'_i-y_i=\alpha_i+\beta_i+1$ for all
$1\leq i\leq r$ since the length of the hook is one plus the sum of
its leg and arm lengths. This relation holds whether or not $X$ is a
minimal $\beta$-set since $y'_i,y_i$ and $\theta(\lambda)$ shift by
the same amount when $X$ is shifted, whereas $\alpha_i$ and
$\beta_i$ do not change. In particular the $\beta$-set $X$ gotten
from removing the largest diagonal hook of $X$ can be used so that
$\theta(\lambda^{\vee})=\theta(\lambda)$ holds. Hence we have the
following.
\begin{lemma} \label{eq_ax} Suppose $\lambda'$ and $\lambda$ are partitions
such that $|\lambda'|<|\lambda|$ and $\lambda'$ can be obtained from
$\lambda$ by removing a sequence of hooks from $\lambda$. Then
$\theta(\lambda)=\theta(\lambda').$
\end{lemma}
\begin{proposition} \label{D-abacus} $X_{-}$ and $X_{+}$ correspond to
$D(\lambda)_{L}$ and
$D(\lambda)_{R}$ in the following manner.  Let $\alpha_i\in
D(\lambda)_{L}$ and $\beta_{i}\in D(\lambda)_{R}$. Then for each
$y'_i\in X_{+}$ and $y_{i}\in X_{-}$ we have
\begin{enumerate}
\item $y_{i}=\theta(\lambda)-\frac{1}{2}-\alpha_i$.
\item $y'_i=\theta(\lambda)+\frac{1}{2}+\beta_i$
\end{enumerate}
\end{proposition}
\begin{proof}
We proceed by induction on $s=|D(\lambda)|$. Suppose $s=1$. Then
$D(\lambda)=\{\alpha|\beta\}$ and $X=\{1,2,\cdots,m-1,m,t_{1}\}$ so
$\theta(\lambda)=m+\frac{1}{2}$. Recall $\alpha$ and $\beta$ are the
number of beads and spaces respectively in the interval $(0,t_{1})$,
and hence $\alpha=\theta(\lambda)-\frac{1}{2}$ and
$\beta=t-\theta(\lambda)-\frac{1}{2}.$

Consider $\lambda$ where $|D(\lambda)|=s$ and let $\lambda^{\vee}$
be the partition obtained by removing $h_{11}$ from $\lambda$. Then,
by induction and Lemma \ref{eq_ax}, when $|D(\lambda^{\vee})|=s-1$
we have that
\[
y_i=\theta(\lambda)-\frac{1}{2}-\alpha_i\;,\;y'_i=\theta(\lambda)+\frac{1}{2}+\beta_i
\]
for $2\leq i\leq r.$ In particular,
$y_2=\theta(\lambda)-\frac{1}{2}-\alpha_2$ and
$y'_2=\theta(\lambda)+\frac{1}{2}+\beta_2.$ But $y_2-y_1$ is one
plus the number of beads between $y_1$ and $y_2$, which is precisely
the difference $\alpha_1-\alpha_2$ by Lemma \ref{bhama}. These
formulas imply $y_1=\theta(\lambda)-\frac{1}{2}-\alpha_1$ and
$y'_1=\theta(\lambda)+\frac{1}{2}+\beta_1$.
\end{proof}
Suppose $\lambda$ is symmetric. Then the diagonal hook lengths
$\delta_{ii}$ are necessarily odd. Then $X$ has a {\it axis of
symmetry} $\theta(\lambda)$ where beads and spaces on one side are
reflected respectively into spaces and beads on the other side.
\begin{corollary} \label{halfint}Suppose $\lambda$ is a symmetric partition
and $X$ is a $\beta$-set in the extended sense for $\lambda$. Then
there exists an axis of symmetry $\theta(\lambda)$ at a half-integer
such that beads and spaces in $X$ to the right of $\theta(\lambda)$
are reflected respectively to spaces and beads in $X$ to the left of
$\theta(\lambda)$.
\end{corollary}
\begin{proof}Follows from Proposition \ref{nsymax} and Proposition
\ref{D-abacus}.
\end{proof}
\begin{lemma}\label{axis} Suppose $\lambda$ is symmetric with empty
$p$-core. Then the number of beads to the right of $\theta(\lambda)$
on the $\gamma$th runner is the same as the number of empty
positions to the left of $\theta(\lambda)$ on the $\gamma$th runner.
\end{lemma}
\begin{proof}
Let $X$ be a $\beta$-set for $\lambda$ and let $X^0$ be the
$\beta$-set obtained from $X$ by sliding all beads to the top of the
$p$-abacus of $X$. Then $X^0$ is a $\beta$-set for $\lambda^0$.
Since $X$ is symmetric about $\theta(\lambda)$, and $X^0$ is
obtained from $X$ by removing successive diagonal hooks, so $X^0$ is
also symmetric about $\theta(\lambda)$. But $X^0$ represents
$\lambda^0=\emptyset$, so it consists of $\{0,1,2,\ldots,t\}$. Then
it is clear that symmetry about $\theta(\lambda)$ is possible only
if $\theta(\lambda)=t+1/2$. Thus all beads of $X^0$ are left of the
axis and all spaces of $X^0$ are right of the axis. Hence, all beads
on the $\gamma$-runner right of the axis must be accommodated by
spaces on the $\gamma$-runner left of the axis.
\end{proof}

We recall that the labeling of the beads on the $\gamma$th runner of
the abacus diagram of $\lambda$ gives a $\beta$-set of the partition
$\lambda_{\gamma}$. Note that each $\lambda_{\gamma}$ has
an axis of symmetry $\theta(\lambda_{\gamma})$. Let $X^0$ be a $\beta$-set of $\lambda^0$
where $|X|=0\pmod{p}.$ In particular, we suppose $|X|=mp$ for some nonnegative integer $m$.

\begin{lemma} \label{mhaf}If $\lambda$ has empty $p$-core, then
$\theta(\lambda_{\gamma})=\theta(\lambda_{\gamma'})$ for all
$0\leq\gamma,\gamma'\leq p-1$. In particular, if $|X|=mp,$ then
$\theta(\lambda_{\gamma})=\theta(\lambda_{\gamma''})=m-\frac{1}{2}.$
\end{lemma}
\begin{proof} Since $\lambda$ has empty $p$-core the $\beta$-set
for $\lambda^0$ is $X=\{0,1,\cdots,mp-1\}$ for some $m$. Then the
abacus diagram will consist of (from north-to-south) $m$ rows of
beads followed by rows of empty spaces. On each runner one begins
counting at 0, hence
$\theta(\lambda^0_{\gamma})=\theta(\lambda^0_{\gamma'})=m-\frac{1}{2}.$
By Lemma \ref{eq_ax},
$\theta(\lambda^0_{\gamma})=\theta(\lambda_{\gamma})$ and
$\theta(\lambda^0_{\gamma})=\theta(\lambda_{\gamma})$. The result
follows.
\end{proof}
We will need the following relation in Section 8.
\begin{lemma} \label{need-shift} Suppose $\lambda$ has empty $p$-core. Then
\[
p(\theta(\lambda_{\gamma})+\frac{1}{2})=\theta(\lambda)+\frac{1}{2}
\]
for all $0\leq \gamma\leq p-1.$
\end{lemma}
\begin{proof} Since $\lambda$ has an empty $p$-core we have $X^0=\{0,1,2,\cdots,mp-1\}$ for some
$m$. Hence $\theta(\lambda)=mp-\frac{1}{2}.$ Since $mp$ is the total
number of beads in $X$, we have
\[
p(m-1+1)=mp-\frac{1}{2}+\frac{1}{2}
\]
which implies
$p(\theta(\lambda_{\gamma})+\frac{1}{2})=\theta(\lambda)+\frac{1}{2}$
by Lemma \ref{mhaf}.
\end{proof}
Let $\lambda$ be such that $\lambda^0\neq \emptyset$ and let
$\bar{\lambda}$ be such that $\bar{\lambda}=\emptyset$ but
$\bar{\lambda}_i=\lambda_i$ for $0\leq i\leq p-1.$
\begin{corollary}\label{adapt} For all $0\leq\gamma\leq p-1$ we have
\[
p(\theta(\bar{\lambda}_{\gamma})+\frac{1}{2})=\theta(\lambda)+\frac{1}{2}.
\]
\end{corollary}
\begin{proof}By Lemma \ref{eq_ax}, $\theta(\bar{\lambda})=\theta(\lambda).$
The result then follows from Lemma \ref{need-shift}.
\end{proof}

We use Corollary \ref{adapt} to offer an interpretation of
$D'(\lambda)$. Suppose $\alpha_i\in D(\lambda)_{L,\gamma}$ so
$\alpha_i=\gamma+\eta_ip.$ Then
$\gamma+\eta_ip=\theta(\lambda)-\frac{1}{2}-y_i$ by Proposition
\ref{D-abacus}.  By Corollary \ref{adapt}
\[
\gamma+\eta_i p = p(\theta(\bar{\lambda}_{\gamma})-\frac{1}{2})+p-1
-y_i.
\]
Suppose $y_i=\gamma^*+{\bar{y}_i}p.$ Then
$\bar{y}_i=\theta(\bar{\lambda}_{\gamma})-\eta_i-\frac{1}{2}.$

Now suppose $\beta_i\in D(\lambda)_{R,\gamma}$ so $\beta_i=\gamma+
\eta_ip$. Then $\gamma+\eta_ip=y'_i-\theta(\lambda)-\frac{1}{2}.$
Then
\[
\gamma+\eta_i p=y'_i-p(\theta(\lambda_{\gamma})+\frac{1}{2}).
\]
by Proposition \ref{D-abacus}. Hence
$y'_i-\gamma=p(\theta(\lambda_{\gamma})+\frac{1}{2}+\eta_i).$
Suppose $y'_i=\gamma+{\bar{y}'_i}p.$ Then
$\bar{y}_i=\theta(\lambda_{\gamma})+\frac{1}{2}+\eta.$ Hence
$D'(\lambda)=\{D_{\gamma}(\lambda)\}_{0\leq \gamma\leq p-1}$ can be
expressed as distances of the beads and spaces of the corresponding
$X_{\gamma}$ from each $\theta(\bar{\lambda}_{\gamma}).$ We will use
this observation in Section 8.

\vskip 15pt
\section{Pairs of straddling or non-straddling $p$-hooks}
If $(x',x]$ is a diagonal hook of $X$ corresponding to a symmetric
$\lambda$, we call $x'$ the {\it opposite position} of $x$. Given
two diagonal hooks $(x',x]$ and $(y',y]$ where $x<y$, we call the
(non-diagonal) hooks $(y',x]$ and $(x',y]$ {\it opposite hooks}.
Conversely, given opposite non-diagonal hooks $(y',x]$ and $(x',y]$
with $x<y$ we get diagonal hooks $(y',y]$ and $(x',x]$.
\begin{lemma} \label{opp res} Suppose $\lambda$ is symmetric and let
$X$ be a $\beta$-set for $\lambda$ with $|X|\equiv 0\pmod{p}$. Let
$(x',x]$ be a diagonal hook of $X$. Then
$x'\equiv{p-1-\gamma}\pmod{p}$ if and only if $x\equiv
\gamma\pmod{p}$.
\end{lemma}
\begin{proof} By symmetry around $\theta(\lambda)$, the number of beads and
empty positions below the axis is $|X|$. Hence
$\theta(\lambda)-\frac{1}{2}\equiv{p-1}\pmod{p}$. Since $x'$ and $x$
are equidistant from $\theta(\lambda)$, then $x'\equiv
p-1-x\pmod{p}$.
\end{proof}
Suppose we want to reduce a symmetric partition $\lambda$ of $n$ to
a symmetric partition $\lambda$ of $n-p$ by removing one $p$-hook.
There is one way of doing so.
\begin{enumerate}
\item ({\it The single hook case}) Then $p$-hook $h=(y,x]$
is a diagonal hook $(x',x]$ where $x-x'=p$.
\end{enumerate}
Suppose we want to reduce a symmetric partition $\lambda$ of $n$ to
a symmetric partition $\lambda'$ of $n-2p$ by removing two
$p$-hooks. By removing two opposite $p$-hooks $h=(y,x]$ and
$h'=(x',y',]$ where $h\neq h'$. There are two cases, the {\it
non-straddling case}, in which $x'<y'<\theta(\lambda)<y<x,$ and the
{\it straddling case}, in which $x'<y<\theta(\lambda)<y'<x$.
\begin{enumerate}
\item {\it (The non-straddling case)}.
Suppose $h$ is completely to the right of $\theta(\lambda)$. Then
removing $h$ and $h'$ is equivalent to replacing a diagonal hook
$(x',x]$ with $(x'+p,x-p]$.

\item {\it (The straddling case)}.
Suppose that $h$ and $h'$ straddle $\theta(\lambda)$. Then removing
$h$ and $h'$ is equivalent to removing two diagonal hooks $(x',x]$
and $(y,y']$ where $x-x'+y'-y=2p$.

\end{enumerate}
Suppose $h=(y,x]$ and $h'=(x',y']$ are non-straddling opposite
$p$-hooks of $\lambda$ (as in Figure 2). Without loss of generality,
$x>y'$. Let $h=h_{ij}$, i.e. have corner $(i,j)$ in $[\lambda]$.
Then the corner of $h$ is on the arm of some diagonal hook. Since
\[
|\{z\in {\mathbb N}:z\not\in X, z\leq y\}|\;>\;|\{z\in {\mathbb N}:
z\in X, z\geq x\}|
\]
by Lemma \ref{bhama}, we have $j>i$. Thus, if $x=\gamma + kp$, where
$0\leq \gamma \leq p-1$ and $k\geq 0$, then $y\not\in X$ such that
$y=\gamma+(k-1)p$. Consequently, $h=(y,x]$ of $\lambda$ corresponds
to some hook $(k-1,k]$ of $\lambda_{\gamma}$ on the $\gamma$th
runner of the $p$-abacus. We give the exact coordinates
$(i_{h},j_{h})$ on the Young diagram $[\lambda_\gamma]$
corresponding to $(k-1,k]$ when $\lambda$ has empty $p$-core. Define
\[
\begin{array}{ccl}
A&=&\{z\equiv\gamma\pmod{p}, z\in X:z \geq x\}\\
B&=&\{z\equiv\gamma\pmod{p},z\not \in X: \theta(\lambda) < z \leq y\}\\
C&=&\{z\equiv -1-\gamma\pmod{p},z\in X: z>\theta(\lambda)     \}.
\end{array}
\]
Let $|A|=a$, $|B|=b$ and $|C|=c$. By construction, $i_h=a$.

By Proposition \ref{nsymax}, $\lambda_{\gamma}$ has an axis
$\theta(\lambda_{\gamma})$ that is a half-integer such that the
number of beads above $\theta(\lambda_{\gamma})$ is the same as the
number of spaces below.
\begin{lemma} \label{noTheta} Suppose $\lambda$ has empty $p$-core and $Y=\{w_1,\cdots,w_j\}$ is
the induced $\beta$-set for $\lambda_{\gamma}$. Let $k$ be an
integer such that $0\leq k\leq w_j$. Then we have the following.
\begin{enumerate}
\item If $k<\theta(\lambda_{\gamma})$ then $(p-\gamma-1)+kp
<\theta(\lambda)$
\item If $k>\theta(\lambda_{\gamma})$ then
$\gamma+kp>\theta(\lambda)$.
\end{enumerate}
\end{lemma}
\begin{proof}
Follows by Proposition \ref{need-shift} and the definition of the
$p$-quotient.
\end{proof}
\begin{lemma}\label{BC} Suppose $\lambda$ is a symmetric partition
with empty $p$-core. Consider the $p$-hook $h=(y,x]$. Let
$(i_h,j_h)$ be the coordinates of the corresponding 1-hook of
$[\lambda_{\gamma}]$ for some fixed $\gamma$. Then
\[
j_h=b\;+\;c
\]
if and only if $h$ is completely to the right of $\theta(\lambda)$.
\end{lemma}
\begin{proof}
Suppose $j_h=b+c$. It is clear that $h$ is completely to the right
of $\theta(\lambda).$ Suppose $h$ is completely to the right of
$\theta(\lambda)$. Since $\lambda$ is symmetric, we know by Lemma
\ref{opp res} that $C$ corresponds bijectively to the set
$\{y'\in{\mathbb N},\; y'\not\in X,\;
y<\theta(\lambda):y'\equiv{\gamma}\pmod{p}\}$. Hence $c$ is also the
number of empty positions less than $\theta(\lambda)$ of residue
$\gamma\pmod{p}$. By Lemma \ref{axis}, $b$ is the number of empty
positions between $\theta(\lambda)$ (and including) $y$ with residue
$\gamma \pmod{p}$. This follows since $\lambda$ has empty $p$-core.
Hence $b+c$ is the total number of empty positions below and
including $y$ with residue $\gamma\pmod{p}$. Then, by Lemma
\ref{bhama}, we are done.
\end{proof}
\begin{lemma} \label{DE} Suppose $\lambda$ is symmetric with empty
$p$-core. Consider the $p$-hook $h=(y,x]$. If $h$ is completely to
the right of $\theta(\lambda)$ then
\[
a\leq c
\]
\end{lemma}
\begin{proof} By Lemma \ref{opp res}, we have that $c$ is the number of empty
positions less than $\theta(\lambda)$ that have residue
$\gamma\pmod{p}$. By Lemma $\ref{axis}$, since $\lambda$ has empty
$p$-core, $c$ must be equal to the number of $z\in X$ such that
$z>\theta(\lambda)$ and $z\equiv\gamma\pmod{p}$. Since
$x>\theta(\lambda)$ and $A=\{z=\gamma + jp,z\in X: z \geq x \}$, we
have $a\leq c$.
\end{proof}
\begin{proposition} \label{arm} A hook $(k-1,k]$ of size 1 on $\lambda_{\gamma}$
corresponds to the $p$-hook $h=(y,x]$ on $\lambda$ where $y$ and $x$
are completely to the right (resp. left) of $\theta(\lambda)$ if and
only if $(k-1,k]$ occurs on an arm (resp. leg) of
$[\lambda_{\gamma}]$.
\end{proposition}
\begin{proof} Suppose $h$ is to the right of
$\theta(\lambda)$. The coordinates of $(k-1,k]$ on the Young diagram
$[\lambda_{\gamma}]$ are $(i_h,j_h)$ where $i_h=a$ (by definition)
and $j_h=b+c$ by Lemma $\ref{BC}$. Since $y\in{B}$, $|B|\neq 0$ and
we have $a<b+c$, since $a\leq c$ by Lemma \ref{DE}. It follows that
$i_h<j_h.$ Hence $(k-1,k]$ is a 1-hook on the arm of
$[\lambda_\gamma]$.

Suppose $(k-1,k]$ is a 1-hook on the arm of $[\lambda_{\gamma}]$.
Clearly $\theta(\lambda_{\gamma})< k-1$. Hence $\theta(\lambda)<
y=\gamma+(k-1)p$ by Lemma \ref{noTheta}. Since $y<x$, $h$ is
completely to the right of $\theta(\lambda)$.

Suppose $(k-1,k]$ is a 1-hook on $\lambda_{\gamma}$ corresponding to
a $p$-hook $h=(y,x]$ completely to the left of $\theta(\lambda)$.
Then $h$ can be viewed as a hook completely to the right of
$\theta(\lambda^*)$. Hence, by the argument above, it corresponds to
a 1-hook on the arm of $\lambda_{\gamma^*}.$  Taking the dual again,
$h$ corresponds to a 1-hook on the leg of $\lambda_{\gamma}.$
\end{proof}
\begin{proposition} \label{str}
Let $\lambda$ be a symmetric partition with empty $p$-core. Let
$\gamma\in [0,p-1]$. Then a hook $h=(k-1,k]$ of size 1 on
$\lambda_{\gamma}$ corresponds to the $p$-hook $h=(y,x]$ of
$\lambda$ where $y$ is to the left of $\theta(\lambda)$ and $x$ is
to the right of $\theta(\lambda)$ if and only if $(k-1,k]$ is a
diagonal hook of $[\lambda_{\gamma}]$.
\end{proposition}
\begin{proof}
Consider a pair $h=(y,x]$ and $h'=(x',y']$ of straddling $p$-hooks.
Without loss of generality, we may suppose $x'<y<y'<x$. Set
$y'=\gamma+jp$ and $x=(p-1-\gamma)+kp$. So $x'=\gamma+(j-1)p$ and
$y=(p-1-\gamma)+(k-1)p$. Suppose $E=\{z\in X,z=\gamma\pmod{p},z\geq
y'\}$ and $F=\{z'\not\in X, z'=\gamma\pmod{p}: z'\leq x'\}$. Let
$|E|=e$ and $|F|=f$. By Lemma \ref{bhama} the coordinates of the
hook $(k-1,k]$ on the Young diagram $[\lambda_{\gamma}]$ are
$(i_h,j_h)=(e,f)$. The inequality $x'<y'<x$ is equivalent to $(k-1)p
< 2\gamma-p+1+jp < kp$. If $\gamma=\frac{p-1}{2}$, then $h=h'$,
which is impossible. Hence we only consider $\gamma\neq
\frac{p-1}{2}$. If $0\le\gamma<(p-1)/2$, it follows that $j=k$ and
$\theta(\lambda)=kp-1/2$. If $(p-1)/2<\gamma\le p-1$, it follows
$j=k-1$ and $\theta(\lambda)=jp+(p-1)/2$. Since $\theta(\lambda)$ is
a half-integer, the second case is impossible and
$0\le\gamma<(p-1)/2$. Now define $E'=\{z\in X:z=\gamma\pmod{p},
z>\theta(\lambda)\}$ and $F'=\{z'\in \mathbb{N}:z'\not\in X,
z'=\gamma\pmod{p}, z'<\theta(\lambda)\}$, so $|E'|=|F'|$ by Lemma
\ref{axis}. Since $y<\theta(\lambda)<x$ and $x-y=p$, we have $\{z\in
X:z=\gamma\pmod{p},\;\theta(\lambda)<z<x\}=\{z'\not\in
X:z'=\gamma\pmod{p},\;y<z'<\theta(\lambda)\}=\emptyset$. Thus $E'=E$
and $F'=F$ and we are done.

Suppose $(k-1,k]$ is a hook of size 1 on the diagonal of
$[\lambda_{\gamma}]$. Then $k-1<\theta(\lambda_{\gamma})<k$. Hence
$(p-1-\gamma+(k-1)p,\gamma+kp]=(y,x]$ straddles $\theta(\lambda)$ by
Lemma \ref{noTheta}.
\end{proof}

\vskip  30 pt
\section{$D'(\lambda)$ concentrated at one or two places}
Suppose $\lambda$ is a symmetric partition. We say $D(\lambda)$ is
{\it concentrated} at $\{\gamma,\gamma^*\}$ if
$D_i(\lambda)\neq\emptyset$ for $i\in \{\gamma,\gamma^*\}$ and
$D_i(\lambda)=\emptyset$ otherwise.
\begin{lemma}\label{dadd} Suppose $\lambda$ and $\lambda'$ are distinct
partitions such that $D(\lambda)$ is concentrated at
$\{\gamma,\gamma^*\}$ and $D(\lambda')$ is concentrated at
$\{\gamma',\gamma'^*\}$ where $\gamma\neq\gamma'$. If
$\{\gamma,\gamma^*\}\neq \{\gamma',\gamma'^*\}$ then $D(\lambda)\cap
D(\lambda')=\emptyset,$ that is, no diagonal hook length of
$\lambda$ equals a diagonal hook length of $\lambda'$.
\end{lemma}
\begin{proof}
Suppose not. Then there exists $\alpha\in D(\lambda)_{L}$ and
$\alpha'\in D(\lambda')_{L}$ such that $\alpha=\alpha'$. But
$\alpha=\gamma+mp$ and $\alpha'=\gamma'+m'p$ so that
$\gamma=\gamma'.$ This is impossible.
\end{proof}
Suppose $\lambda$ and $\lambda'$ are symmetric partitions such that
$D(\lambda)\cap D(\lambda')=\emptyset.$  Define $\lambda+\lambda'$
to be the symmetric partition such that
$D(\lambda+\lambda')=D(\lambda)\cup D(\lambda').$ In particular, we
can form $\lambda+\lambda'$ whenever $\lambda$ and $\lambda'$ are
concentrated on disjoint sets.
\begin{theorem}\label{eqlgth} Let $\lambda$ be symmetric with empty
$p$-core such that $D'(\lambda)$ is concentrated at
$\{\gamma,\gamma^*\}$ where $\gamma\neq\gamma^*$. Then
\begin{enumerate}
\item $D'(\lambda)$ is a $p$-tuple of bisequences, that is, for each
$\gamma$, $D_{\gamma}(\lambda)_{L}$ and $D_{\gamma}(\lambda)_{R}$
are of equal lengths.
\item For each $\gamma$, $D_{\gamma}(\lambda)=D(\lambda_{\gamma})$
and $D_{\gamma^*}(\lambda)=D(\lambda_{\gamma^*}),$ where
$\lambda_{\gamma}$ and $\lambda_{\gamma^*}$ are the $\gamma$th and
$\gamma^*$th components of the $p$-quotient of $\lambda$.
\item Suppose
$D(\lambda_{\gamma})=(\sigma_{1},\cdots,\sigma_{w}|\tau_{1},\cdots,\tau_{w})$.
Then
\[
D(\lambda)=(\alpha_{1},\cdots,\alpha_{2w}|\alpha_{1},\cdots,\alpha_{2w})
\]
where $\{\alpha_1,\cdots,\alpha_{2w}\}=$
\[
\{\gamma^*+\sigma_{i}p \;, \gamma+\tau_{i}p:\;1\leq i\leq w\}.
\]
\end{enumerate}
\end{theorem}
\begin{proof}
By induction on $|\lambda|$. The minimal case is $|\lambda|=2p$
where $D(\lambda)=(p-1-\gamma,\gamma|p-1-\gamma,\gamma)$. Then
$\lambda$ is comprised of just two opposite $p$-hooks. Hence
$|D_{\gamma}(\lambda)_{R}|=|D_{\gamma}(\lambda)_{L}|=1$ and part (1)
follows. By definition, $D_{\gamma}(\lambda)=(0|0)$ and
$D_{p-1-\gamma}(\lambda)=(0|0)$. Since $\lambda_{\gamma}=(1)$ and
$\lambda_{\gamma^*}=(1)$, part (2) follows. Part (3) follows since
$D(\lambda)=(p-1-\gamma,\gamma|p-1-\gamma,\gamma)$. Now suppose
$|\lambda|=n>2p.$ By induction, we assume that the theorem holds for
all partitions $\lambda$ such that $|\lambda|<n$. Consider
$|\lambda|=n$. Let $h,h'$ be opposite $p$-hooks in $\lambda$ and let
$\lambda^{\vee}$ be the symmetric partition gotten from removing $h$
and $h'$. Following the discussion preceding Lemma \ref{BC}, there
are two cases.\\Case 1: (The non-straddling case) Here one obtains
$D(\lambda^{\vee})$ from $D(\lambda)$ by replacing an element
$(\alpha|\alpha)$ by $(\alpha-p|\alpha-p)$ where $\alpha-p\geq 0$.
Then $D_{\gamma}(\lambda)_L$ and $D_{\gamma}(\lambda^{\vee})_L$ are
the same except for some $\sigma_{\mu}\in D_{\gamma}(\lambda)_L$
which is replaced by $\sigma_{\mu}-1$. By symmetry, $\sigma_{\mu}$
is replaced by $\sigma_{\mu}-1$ resulting in
$D_{\gamma^*}(\lambda')_R$. We prove that parts (1), (2), and (3)
hold.
\begin{enumerate}
\item  By induction $D_{\gamma}(\lambda^{\vee})$ and
$D_{\gamma^*}(\lambda^{\vee})$ are both bisequences with components
of equal length. Hence the same is true for $D_{\gamma}(\lambda)$
and $D_{\gamma^*}(\lambda)$.
\item By induction
$D_{\gamma}(\lambda^{\vee})=D(\lambda^{\vee}_{\gamma})$. By
Proposition \ref{arm}$, \lambda_{\gamma}$ is obtained from
$\lambda^{\vee}_{\gamma}$ by adding a hook of size 1 to both the leg
length of the diagonal hook corresponding to
$(\sigma_{\mu}-1|\tau_{\mu})\in D(\lambda'_{\gamma})$ and to the arm
length of the diagonal hook $(\tau_{\mu}|\sigma_{\mu}-1)\in
D(\lambda^{\vee}_{\gamma^*})$. Hence
$D_{\gamma}(\lambda)=D(\lambda_{\gamma})$ and
$D_{\gamma^*}(\lambda)=D(\lambda_{\gamma^*})$.
\item Given
\[
D(\lambda^{\vee}_{\gamma})=(\sigma_{1},\cdots,\sigma_{\mu}-1,\cdots,\sigma_{w}|\tau_{1},\cdots,\tau_{i},\cdots,\tau_{w})
\]
\[
D(\lambda^{\vee}_{p-1-\gamma})=(\tau_{1},\cdots,\tau_{i},\cdots,\tau_{w}|\sigma_{1},\cdots,\sigma_{\mu}-1,\cdots,\sigma_{w})
\]
we have by induction that
\[
D(\lambda^{\vee})=(\cdots\alpha'_{i},\alpha''_{i}\cdots|\cdots\alpha'_{i},\alpha''_{i}\cdots)
\]
where
\[
\begin{array}{cc}
\alpha'_{i}=(p-1-\gamma)+\sigma_{i}p\qquad & \qquad\alpha''_{i}=\gamma+\tau_{i}p\\
\end{array}
\]
$1\leq i \leq w$ and $i\neq\mu$. When $i=\mu$, then
\[
\begin{array}{cc}
\alpha'_{\mu}=(p-1-\gamma)+(\sigma_{\mu}-1)p\qquad & \qquad\alpha''_{\mu}=\gamma+\tau_{\mu}p\\
\end{array}
\]
It is clear by replacing $\sigma_{\mu}-1$ by $\sigma_{\mu}$, that
the desired formula for $D(\lambda)$ is obtained.
\end{enumerate}
Case 2:(The straddling case) Here $D(\lambda^{\vee})$ one obtains
$D(\lambda)$ from by removing $(\alpha|\alpha)$ and $(\beta|\beta)$
where $\alpha+\beta+1=p$ (assume without loss of generality that
$\alpha>\beta$). Then $D(\lambda^{\vee})$ is also concentrated at
$\{\gamma,\gamma^*\}$. The relation $\alpha+\beta+1=p$ implies that
if $\alpha=\gamma$, then $\beta=p-1-\gamma$. Thus $(\alpha|\alpha)$
contributes a term 0 to $D_{\gamma}(\lambda)_{L}$ and a 0 to
$D_{\gamma^*}(\lambda)_{R}$. Likewise $(\beta|\beta)$ contributes a
term 0 to $D_{\gamma}(\lambda)_R$ and a 0 to
$D_{\gamma^*}(\lambda)_{R}$. We prove that parts (1), (2), and (3)
hold.
\begin{enumerate}
\item By induction $D_{\gamma}(\lambda^{\vee})$ is a bisequence with
components of equal length. Hence $D_{\gamma}(\lambda)$.
\item By
induction $D_{\gamma}(\lambda^{\vee})=D(\lambda^{\vee}_{\gamma})$
and $D_{p-1-\gamma}(\lambda^{\vee})=D(\lambda^{\vee}_{\gamma^*})$.
Now we re-attach to $\lambda^{\vee}$ the diagonal hooks
corresponding to $(\alpha|\alpha)$ and $(\beta|\beta)$, where
$\alpha+\beta+1=p.$ This is equivalent to adjoining $(0|0)$ to both
$D_{\gamma}(\lambda)$ and $D_{\gamma^*}(\lambda')$. The effect on
the partitions $\lambda^{\vee}_{\gamma}$ and
$\lambda^{\vee}_{\gamma^*}$ will be adding a diagonal node of size 1
to each, by Proposition \ref{str}. Hence
$D_{\gamma}(\lambda)=D(\lambda_{\gamma})$ and
$D_{\gamma^*}(\lambda)=D(\lambda_{\gamma^*})$.
\item Given
\[
\begin{array}{ccc}
D(\lambda_{\gamma}^{\vee})&=&(\sigma_{1},\cdots,\sigma_{w-1}|\tau_{1},\cdots,\tau_{w-1})\\
D(\lambda_{\gamma^*}^{\vee})&=&(\tau_{1},\cdots,\tau_{w-1}|\sigma_{1},\cdots,\sigma_{w-1})
\end{array}
\]
then by induction
\[
D(\lambda^{\vee})=(\cdots\alpha'_{i},\alpha''_{i}\cdots|\cdots\alpha'_{i},\alpha''_{i}\cdots)
\]
where
\[
\begin{array}{cc}
\alpha'_{i}=(p-1-\gamma)+\sigma_{i}p\qquad & \qquad\alpha''_{i}=\gamma+\tau_{i}p\\
\end{array}
\]
and $1\leq i\leq w-1$. Now attaching $(0|0)$ to
$D(\lambda^{\vee}_{\gamma})$ and $(0|0)$ to
$D(\lambda^{\vee}_{\gamma^*})$, is equivalent to adjoining both
$p-1-\gamma$ and $\gamma$ to both $D(\lambda)_R$ and $D(\lambda)_L$.
Hence the desired formula for $D(\lambda)$ is obtained.
\end{enumerate}
\end{proof}
\begin{corollary} \label{Cor1} Suppose $\lambda$ is symmetric with empty $p$-core and $D(\lambda)$ is concentrated at
$\{\gamma,\gamma^*:\gamma\neq\gamma^*\}$ and $D_{\gamma}(\lambda)=(\sigma_1,\cdots,\sigma_w|\tau_1,\cdots,\tau_w)$. Then
\[
\delta(\lambda)=\cup_{i}\{2(\sigma_i+1)p-2\gamma-1,2{\tau_i}p+2\gamma+1\}
\]
\end{corollary}
\begin{proof} Follows from part 3 of Theorem \ref{eqlgth}.
\end{proof}
\begin{example} Let $p=5$. Suppose the $5$-quotient is concentrated at $\{\gamma,\gamma^*\}=\{0,4\}$, where
$\lambda_{0}=(6^2,2)$ and $\lambda_{4}=(3^2,2^3)$. Then
$D_0(\lambda)=(2,1|5,4)$ and $D_4(\lambda)=(5,4|2,1)$,
$D(\lambda)=(25,20,14,9|25,20,14,9)$ and
$\delta(\lambda)=(51,41,29,19)$.
\end{example}

Similar results hold for the case when $\lambda$ is symmetric and
$D(\lambda)$ is concentrated at $\gamma=\gamma^*=\frac{p-1}{2}$.
\begin{theorem}\label{oneplace} Suppose $\lambda$ is a symmetric partition
with empty $p$-core and let $D(\lambda)$ be concentrated at
$\gamma=\gamma^*=\frac{p-1}{2}$. Then
\begin{enumerate}
\item $D'(\lambda)$ is a $p$-tuple of $p-1$ empty bisequences, with
$D_{\frac{p-1}{2}}(\lambda)\neq \emptyset$ and
${D_{\frac{p-1}{2}}}(\lambda)_R$ and
${D_{\frac{p-1}{2}}}(\lambda)_L$ are of equal lengths.
\item $D_{\frac{p-1}{2}}(\lambda)=D(\lambda_{\frac{p-1}{2}})$.
\item Suppose
$D(\lambda_{\frac{p-1}{2}})=(w_{1},\cdots,w_{\mu}|w_{1},\cdots,w_{\mu})$,
and $D(\lambda_{\gamma})=\emptyset$ when $\gamma\neq\frac{p-1}{2}$.
Then \[ D(\lambda)=(z_1,\cdots,z_{\mu}|z_{1},\cdots,z_{\mu}) \]
where $z_{i}=\frac{p-1}{2}+w_{i}p$.
\end{enumerate}
\end{theorem}
\begin{proof}
By induction on $|\lambda|$. The minimal case is $|\lambda|=p$. In
this case $D(\lambda)=(\frac{p-1}{2}|\frac{p-1}{2})$ and
$D_{\frac{p-1}{2}}(\lambda)=(0|0)$. The remainder of the proof is
similar to that of Theorem $\ref{eqlgth}$.
\end{proof}
\begin{corollary} \label{oneplc2} Suppose $\lambda$ is symmetric with empty $p$-core, such
that $\lambda$ is concentrated at $\{\frac{p-1}{2}\}$. Then
$\delta(\lambda)=\cup_{i}\{(2m_i+1)p\}$ if
$\delta(\lambda_{\frac{p-1}{2}})=\cup_{i}\{2m_i+1\}$ for every
$(m_i|m_i)\in D(\lambda_{\frac{p-1}{2}})$.
\end{corollary}
\begin{proof}This follows from Theorem 2, part 3.
\end{proof}
\begin{example}Let $p=5$. Suppose the $5$-quotient is concentrated at
$\{2\}$ and $\lambda_{2}=(2^2)$. Then $D_2(\lambda)=(1,0|1,0)$ and
$D(\lambda)=(7,2|7,2)$, $\delta(\lambda)=(15,5)$.
\end{example}

{\vskip 30pt}
\section{Symmetric partitions with an empty $p$-core}
Now suppose $\lambda$ is symmetric and has empty $p$-core. Fix a
$\gamma$ between $0$ and $\frac{p-1}{2}$. Suppose
$D(\lambda_{[\gamma]})\subseteq D(\lambda)$ is the bisequence whose
$p$-quotient $D'(\lambda_{[\gamma]})$ has just the components
$D_{\gamma}(\lambda)$ and $D_{\gamma^*}(\lambda)$. Let
$\lambda_{[\gamma]}$ be the symmetric partition corresponding to
$D(\lambda_{[\gamma]})$. By Lemma \ref{dadd},
\[
\lambda_{[0]},\cdots,\lambda_{[\frac{p-1}{2}]}
\]
have disjoint diagonals. Thus
$\lambda_{[0]}+\lambda_{[1]}+\cdots+\lambda_{[\frac{p-1}{2}]}$ is
defined in the sense described in the remark before Theorem
\ref{eqlgth}.
\begin{theorem}\label{main4} Suppose $\lambda$ is symmetric
and has empty $p$-core.
Then
\[
\lambda=\lambda_{[0]}+\lambda_{[1]}+\cdots+\lambda_{[\frac{p-1}{2}]}
\]
\[
D(\lambda)=\coprod_{1\leq\gamma\leq\frac{p-1}{2}}D(\lambda_{[\gamma
]})
\]
\end{theorem}
\begin{proof}
By Lemma \ref{dadd}, it is clear that $D(\lambda_{[\gamma]})\cap
D(\lambda_{[\mu]})=\emptyset$ when $\gamma\neq\mu$. Let
$k_{\gamma}=|D_{\gamma}(\lambda)_R|$ By Theorem \ref{eqlgth} and
Theorem \ref{oneplace}, for a fixed $\gamma$, we have for all $1\leq
i\leq k_{\gamma}$ and $1\leq j\leq t$, the diagonal hooks of
$\lambda$ corresponding to
$(\alpha'_{\gamma,i}|\alpha'_{\gamma,i})$, $(\alpha''_{\gamma,i}|
\alpha''_{\gamma,i})$ and $(z_j|z_j)$ in $D(\lambda)$ have distinct
lengths. Hence
$\bigcap_{1\leq\gamma\leq\frac{p-1}{2}}D(\lambda_{[\gamma]})=\emptyset$.
Since these exhaust the diagonal hooks arising from the $p$-quotient
$D'(\lambda)$, and $\lambda$ has an empty $p$-core,
$\coprod_{1\leq\gamma\leq\frac{p-1}{2}}D(\lambda_{[\gamma ]})$
constitute all of the diagonal hook lengths of $\lambda$.
\end{proof}
\begin{example} \label{empty general example}
\end{example}
Suppose $p=5$, $\lambda\vdash 190$ is symmetric with empty $p$-core
and $\lambda_{0}=(6^2,2), \lambda_1=(3), \lambda_2=(2^2),
\lambda_3=(1^3), \lambda_4=(3^2,2^3)$. Then
$D_{0}(\lambda)=\{2,1|5,4\}, D_{1}(\lambda)=\{0|2\},
D_{2}(\lambda)=\{1,0|1,0\}, D_{3}(\lambda)=\{2|0\}$, and
$D_{4}(\lambda)=\{5,4|2,1\}$. Hence
\[
D(\lambda)=\{25,20,14,11,9,7,3,2|25,20,14,11,9,7,3,2\}
\]
and $\delta(\lambda)=(51,41,29,23,19,15,7,5)$.
\section{Symmetric $p$-cores}
For any partition $\lambda$, let
\[ D(\lambda)_{L,\gamma}=\{\alpha\in
D(\lambda)_L:\alpha\equiv\gamma\pmod{p}\}
\]
and
\[
D(\lambda)_{R,\gamma}=\{\beta\in
D(\lambda)_R:\beta\equiv\gamma\pmod{p}\}.
\]
Let $\lambda^0$ be a symmetric $p$-core partition. Let
$D(\lambda)_{\gamma}$ be the set of $(\beta|\beta)\in D(\lambda)$
such that $\beta\equiv \gamma\pmod{p}$. Then, in particular,
$D(\lambda^0)=\cup_{\gamma}D(\lambda^0)_{\gamma}$ and
$D(\lambda^0)_{\gamma'}\cap D(\lambda^0)_{\gamma}=\emptyset$ for
$\gamma\neq\gamma'$.
\begin{proposition} \label{empty_star}Suppose $\lambda^0$ is a symmetric $p$-core
partition. Then for $\gamma\neq\gamma^*$,
\[
D(\lambda^0)_{\gamma}\neq\emptyset\;\text{ implies }
D(\lambda^0)_{\gamma^*}=\emptyset.
\]
\end{proposition}
\begin{proof} Suppose $(\alpha_i|\alpha_i)\in D(\lambda^0)_{\gamma}$
and $(\beta_j|\beta_j)\in D(\lambda^0)_{\gamma^*}$. Then in the
notation of Proposition \ref{D-abacus} we have $y'_i\equiv
\theta(\lambda^0)+\frac{1}{2}+\gamma\pmod{p}$ and
$y_j\equiv\theta(\lambda^0)-\frac{1}{2}-\gamma^*\pmod{p}.$ But
$\gamma^*=-1-\gamma\pmod{p}$. Thus $y'_i-y_j\equiv\pmod{p}$,
contradicting the assumption $\lambda^0$ is a $p$-core.
\end{proof}
A symmetric $\lambda$ is {\it $\gamma$-packed} if
$D(\lambda)_{\gamma}$ consists of the elements
$(\gamma+ip|\gamma+ip)$ for $i=0,1,\cdots,r$. Let $X_{\gamma, +}$ be
the subset of $X_+$ consisting of elements $y'$ where
\[
y'-\theta(\lambda)-\frac{1}{2}\equiv \gamma \pmod{p}.
\]
We define $X_{\gamma,-}$ similarly.
\begin{proposition} \label{packed} Suppose $\lambda^0$ is a symmetric $p$-core and
$D(\lambda)_{\gamma}\neq\emptyset$. Then $\lambda^0$ is
$\gamma$-packed.
\end{proposition}
\begin{proof} Clearly if $\lambda^0$ is not $\gamma$-packed, then
there exist integers $y'$,$z'$, greater than $\theta(\lambda^0)$
such that $z'\equiv y'\pmod{p}$, $z'\not\in X_+$, and $y'\in X_+$.
In particular, $(z',y']$ is a $p$-hook of $\lambda^0$.
\end{proof}
\begin{corollary}(symmetric $p$-core criterion)
\label{p-core criterion}Let $\lambda^0$ be a symmetric partition.
Then $\lambda^0$ is a $p$-core if and only if for every
$\gamma\in\{0,\cdots,p-1\}$ $D(\lambda^0)_{\gamma}\neq\emptyset$
implies that $\lambda^0$ is $\gamma$-packed and
$D(\lambda^0)_{\gamma^*}=\emptyset$.
\end{corollary}
\begin{proof} Clearly, if $\lambda^0$ if a $p$-core the result follows
by Proposition \ref{empty_star} and Proposition \ref{packed}.
Suppose that for each $\gamma\in\{0,\cdots,p-1\}$, if
$D(\lambda^0)_{\gamma}\neq\emptyset$ then $\lambda^0$ is
$\gamma$-packed and $D(\lambda^0)_{\gamma^*}=\emptyset$, but that
$\lambda^0$ is not a $p$-core. Then, for some $\gamma$ there exists
a hook $h=(x,y']$ where $y'=\gamma+mp$, $x=\gamma+(m-1)p$. By
symmetry, we can assume that $y'>\theta(\lambda^0)$. If
$x>\theta(\lambda)$, then $\lambda$ is not $\gamma$-packed, which is
a contradiction. Now suppose $x<\theta(\lambda^0)$, then by symmetry
there exists $x'>\theta(\lambda^0)$ such that
$x'-\theta(\lambda^0)-\frac{1}{2}\equiv \gamma^*\pmod{p}$. This
implies $D(\lambda^0)_{\gamma^*}\neq\emptyset$, which is a
contradiction.
\end{proof}
\begin{example} \label{p-core general example}
\end{example}
Suppose $p=5$ and $\lambda'\vdash 324$ such that $\lambda'$ is
symmetric and
\[
\delta(\lambda')=(69,59,49,39,29,27,19,17,9,7).
\]
In particular, $D(\lambda')_{R}=(34,29,24,19,14,13,9,8,4,3)$. Hence
$D(\lambda')_{R,4}=(34,29,24,19,14,9,4)$ and
$D(\lambda')_{R,3}=(13,8,3)$. Hence $\lambda'$ is both 4-packed and
3-packed. Since $D(\lambda')_0=\emptyset$ and
$D(\lambda')_1=\emptyset$, $\lambda'$ is a $5$-core by Theorem
\ref{p-core criterion}.

{\vskip 30pt}
\section{Symmetric partitions with a non-empty $p$-core}
We extend the results of Section 6 to the case of a symmetric
partition with a non-empty $p$-core. Let $\bar{\lambda}$ be the
symmetric partition that shares the $p$-quotient with $\lambda$, but
has empty $p$-core. Hence $(\bar{\lambda})^{0}=\emptyset$ and
$(\bar{\lambda})_{\gamma}=\lambda_{\gamma}$ for $0\leq \gamma\leq
p-1.$

Now consider a symmetric partition $\lambda$ of $n$ with a non-empty
$p$-core $\lambda^0$. Let $\bar{X}$ and $X^0$ be $\beta$-sets of
$\bar{\lambda}$ and $\lambda^0$ respectively. Since
$\lambda^0\neq\emptyset$, we have $(X^0)_{\gamma,+}\neq\emptyset$ ,
for some $\gamma$. Then $|D(\lambda^0)_{\gamma}|\neq\emptyset$. In
particular, $|D(\lambda^0)_{\gamma}|=d^{0}_{\gamma}$ by Proposition
\ref{empty_star}. By the definition of $D'(\lambda)$ each
$(\alpha|\alpha)\in D(\lambda^0)_{\gamma}$ contributes an element to
both $D_{\gamma^*}(\lambda)_L$ and $D_{\gamma}(\lambda)_R$.
($D(\lambda^0)$ contributes nothing to $D_{\gamma^*}(\lambda)_R$
and $D_{\gamma}(\lambda)_L$.) The definition of $D'(\lambda)$ (and
Proposition \ref{empty_star}) forces
$|D_{\gamma}(\lambda)_R|-|D_{\gamma}(\lambda)_L|=d^{0}_{\gamma}$.
This implies $D'(\lambda)$ is not a $p$-tuple of bisequences.
Specifically, $D_{\gamma}(\lambda)\neq D(\lambda_{\gamma})$.

Define ${\Omega'}\subset\{0,\cdots,p-1\}$ so that $\gamma'\in
{\Omega'}$ if $D_{\gamma'}(\lambda^0)\neq\emptyset$ (i.e.
$d^0_{\gamma'}> 0$).  Let $(\Omega')^*=\{p-\gamma'-1:\gamma'\in \Omega'\}$
and $U=\Omega'\cup(\Omega')^*.$  Define $\Omega''=\{0,\cdots,p-1\}-U.$

\begin{lemma}\label{pushdown}
\noindent
\begin{enumerate}
\item $\theta(\lambda_{\gamma''})=\theta(\bar{\lambda}_{\gamma''})$
\item $\theta(\lambda_{\gamma'})=\theta(\bar{\lambda}_{\gamma'})+d^0_{\gamma'}$
\end{enumerate}
\end{lemma}
\begin{proof}
When all beads on all the runners of the abacus of $\lambda$ are
moved up completely one obtains the abacus diagram for $\lambda^0$
(Theorem 2.7.16, \cite{J-K}). Since $D(\lambda^0)_{\frac{p-1}{2}}$
is empty, $d^0_{\frac{p-1}{2}}=0$ and the $\frac{p-1}{2}$th runner
of $\lambda^0$ is unchanged from the $\frac{p-1}{2}$th runner of
$\bar{\lambda}^0$. Let $\bar{X}^0$ and $X^0$ be the $\beta$-sets for
$\lambda^0$ and $\bar{\lambda}^0.$ Matching the $\frac{p-1}{2}$th
runners of $\lambda^0$ and $\bar{\lambda}^0$ one can superimpose the
abacus of $\lambda^0$ onto of the abacus of $\bar{\lambda}^0$. It
follows that $|\bar{X}^0_{\gamma'}|+d^0_{\gamma}=|X^0_{\gamma'}|$
for $\gamma'\in\Omega'.$ Also,
$|\bar{X}^0_{\gamma''}|=|X^0_{\gamma''}|$ since $d^{0}_{\gamma''}=0$
for $\gamma''\in \Omega''$. Hence
$\theta(\lambda^0_{\gamma'})=\theta(\bar{\lambda}^0_{\gamma'})+d^0_{\gamma'}$
and
$\theta(\lambda^0_{\gamma''})=\theta(\bar{\lambda}^0_{\gamma''})$.
The result follows since
$\theta(\lambda_{\gamma'}^0)=\theta(\bar{\lambda}_{\gamma'}^0)$ and
$\theta(\lambda_{\gamma''}^0)=\theta(\bar{\lambda}_{\gamma''}^0)$ by
Proposition \ref{eq_ax}.
\end{proof}

We can describe $X_{\gamma'}$ using $\bar{X}_{\gamma'}$ and Lemma
\ref{pushdown} in the following three steps which we call the {\it
$d^0_{\gamma'}$-shift} of $\bar{X}_{\gamma'}.$
\begin{enumerate}
\item $m_{\sigma}\in X_{\gamma',+}$ if $m_{\sigma}-d^0_{\gamma'}>\theta(\bar{\lambda}_{\gamma})$
and $m_{\sigma}-d^0_{\gamma'}\in \bar{X}_{\gamma',+}$\\
\item $m_s\in X_{\gamma',+}$ if $\theta(\bar{\lambda}_{\gamma'})<m_s<\theta(\bar{\lambda}_{\gamma'})+d^0_{\gamma'}$
and $m_s-d^0_{\gamma'}\not\in \bar{X}_{\gamma',-}$\\
\item $m_t\in X_{\gamma',-}$ if
$m_t-d^0_{\gamma'}<\theta(\bar{\lambda}_{\gamma'})$.
\end{enumerate}

Now consider the following sets
\[
\begin{array}{lll}
S_{\gamma'}(\bar{\lambda})_L&=:&\{s:s\in{\mathbb N}, s\not\in
D_{\gamma'}(\bar{\lambda})_{L},0\leq s\leq d^0_{\gamma'}-1\}\\
T_{\gamma'}(\bar{\lambda})_L&=:&\{t: t\in
D_{\gamma'}(\bar{\lambda})_L, t\geq d^{0}_{\gamma'}\}.
\end{array}.
\]
Following the comments after Proposition \ref{adapt},
$S_{\gamma'}(\bar{\lambda})_L$ and $T_{\gamma'}(\bar{\lambda})_L$
are in bijection with the subsets of $\bar{X}_{\gamma}$ in steps (2)
and (3) of the definition of the $d^0_{\gamma}$-shift.  Now we can
interpret $D_{\gamma'}(\lambda)$ via the $d^0_{\gamma'}$-shift of
$\bar{X}_{\gamma'}.$
\begin{proposition} \label{DD-shift} $D_{\gamma'}(\lambda)$ is obtained from
$D_{\gamma'}(\bar{\lambda})$ in the following three steps.
\begin{enumerate}
\item Each  $\sigma\in
D_{\gamma'}(\bar{\lambda})_R$ is sent to
$\sigma+d^0_{\gamma}\in D_{\gamma'}(\lambda)_R$
\item Each $s\in S_{\gamma'}(\bar{\lambda})_L$ is sent to $d^0_{\gamma'}-s-1\in D_{\gamma'}(\lambda)_R$
\item Each $t\in T_{\gamma'}(\bar{\lambda})_L$ is sent to $t-d^0_{\gamma'}\in
D_{\gamma'}(\lambda)_L$.
\end{enumerate}
\end{proposition}
\begin{proof}
\noindent We prove part (2). Let $x_s=m_s-d^0_{\gamma'}$. Each
$x_s<\theta(\bar{\lambda}_{\gamma})$ where $x_s\not\in
\bar{X}_{\gamma,-}$ where
$x_s+d^0_{\gamma}>\theta({\bar{\lambda}_{\gamma}})$ corresponds to
some $s\in S_{\gamma'}(\lambda)_L$. Hence we have
$\gamma+(x_s+d^0_{\gamma})p\in X_+$ by the usual $p$-quotient.
Again, by Proposition \ref{D-abacus},
\[
\theta(\lambda)+\frac{1}{2}+\beta=\gamma+(x_s+d^0_{\gamma'})p
\]
for some $\beta\in  D(\lambda)_R.$ By substitution,
\[
\beta=\gamma-\theta(\lambda)-\frac{1}{2}+(\theta(\bar{\lambda}_{\gamma})-\frac{1}{2})p+(d^0_{\gamma'}-s)p.
\]
By Lemma \ref{need-shift} we have
\[
\beta=\gamma+(d^{0}_{\gamma'}-s-1)p.
\]
By definition of $D'(\lambda)$, $x_s+d^0_{\gamma}$ corresponds to
$d^{0}_{\gamma'}-\bar{s}-1\in D_{\gamma'}(\bar{\lambda})_R$.

The proofs of (1) and (3) are similar.
\end{proof}
\begin{proposition} \label{nonempty_1}
Given $\lambda^0$ and $\bar{\lambda}$,
\[
D(\lambda)_{R}=O_1\cup O_2\cup O_3 \cup O_4
\] where
\[
\begin{array}{ccll}
O_1&=&\cup_{\gamma'\in \Omega'}
&\{\gamma'+(\sigma+d^0_{\gamma'})p:\sigma\in
D_{\gamma'}(\bar{\lambda})_{R}\}\\
O_2&=&\cup_{\gamma'\in \Omega'}
&\{\gamma'+(d^0_{\gamma'}-s-1)p:s\in S_{\gamma'}(\bar{\lambda})\}\\
O_3&=&\cup_{\gamma'\in \Omega'}
&\{(p-1-\gamma')+(t-d^0_{\gamma'})p:t\in T_{\gamma'}(\bar{\lambda})\}\\
O_4&=&\cup_{\gamma''\in \Omega''} &\{\gamma'+\mu p:\mu\in
D_{\gamma''}(\bar{\lambda})_{R}\}.
\end{array}
\]
\end{proposition}
\begin{proof} This follows from the definition of $D'(\lambda)$ and
Proposition \ref{DD-shift}.
\end{proof}
\begin{theorem} \label{nonempty_2}
Given $\lambda^0$ and $\bar{\lambda}$,
\[\delta(\lambda)={\mathbb
O}_1\cup{ \mathbb O}_2\cup{\mathbb O}_3\cup {\mathbb O}_4 \] where
\[
\begin{array}{ccll}
{\mathbb O}_1&=&\cup_{\gamma'\in
{\Omega'}}&\{2(\gamma'+(\sigma+d^0_{\gamma'})p)+1:\sigma\in
D_{\gamma'}(\bar{\lambda})_{R}\}\\
{\mathbb O}_2&=&\cup_{\gamma'\in \Omega'}
&\{2(\gamma+(d^0_{\gamma'}-s-1)p)+1:s\in S_{\gamma'}(\lambda)_L\}\\
{\mathbb O}_3&=&\cup_{\gamma'\in \Omega'}
&\{2((p-1-\gamma')+(t-d^0_{\gamma'})p)+1:t\in
T_{\gamma'}(\lambda)_L\}.\\
{\mathbb O}_4&=&\cup_{\gamma''\in\Omega''} &\{2(\gamma'+\mu
p)+1:\mu\in D_{\gamma''}(\bar{\lambda})_{R}\}
\end{array}
\]
where $d^0_{\gamma'}=|D(\lambda^0)_{\gamma'}|,$ $\gamma'\in \Omega'$
and $D_{\gamma'}(\lambda)_R,$ $S_{\gamma'}(\lambda)_L$ and
$T_{\gamma'}(\lambda)_L$ are as above.
\end{theorem}
\begin{proof} Follows from the Proposition \ref{nonempty_1} and the relationship between
$D(\lambda)$ and $\delta(\lambda)$.
\end{proof}
[Note: In the case $d^0_{\gamma}=0$ for all $0\leq \gamma\leq p-1$,
Proposition \ref{nonempty_1} reverts to Theorem \ref{main4}.]
\begin{example} \label{nonempty general example}
\end{example}
Suppose $p=5$ and $\eta\vdash 514$ such that $\eta$ is symmetric
such that $\eta^0=\lambda'$, where $\lambda'$ is as in Example
$\ref{p-core general example}$. Furthermore, let
$\eta_{i}=\lambda_i$ for $0\leq i\leq p-1$, where $\lambda_i$ is as
in Example \ref{empty general example}. In this case we have
$d^0_3=3,$ $d^0_4=7,$ $D_{3}(\bar{\eta})=\{2|0\},$ and
$D_{4}(\bar{\eta})=\{5,4|2,1\}$. Hence, by
 Proposition
\ref{DD-shift}, we have
\[
\begin{array}{lcc}
D_{2}(\eta)_{R}&=&\{1,0\}\\
D_{2}(\eta)_{L}&=&\{1,0\}\\
D_{3}(\eta)_{R}&=&\{2,1,0\}\\
D_{3}(\eta)_{L}&=&\emptyset\\
D_{4}(\eta)_{R}&=&\{9,8,6,5,4,3,0\}\\
D_{4}(\eta)_{L}&=&\emptyset.\\
\end{array}
\]
Then by Proposition \ref{nonempty_1},
\[
\begin{array}{lcc}
D(\eta)_{R,2}&=&\{6,2\}\\
D(\eta)_{R,3}&=&\{13,8,3\}\\
D(\eta)_{L,3}&=&\emptyset\\
D(\eta)_{R,4}&=&\{49,44,34,29,24,19,4\}\\
D(\eta)_{L,4}&=&\emptyset.\\
\end{array}
\]
Finally, by Theorem \ref{nonempty_2} we have
\[
\delta(\eta)=(99,69,59,49,39,37,27,17,13,9,7,5).
\]
{\bf Acknowledgements} Portions of this paper appeared in my PhD thesis \cite{Nat_thesis} completed under the direction
 of Paul Fong.  The author would also like to thank Itaru Terada for helpful discussions.

\end{document}